\documentclass[12pt]{article}
\usepackage{amsmath}
\usepackage{color}
\usepackage{fancyhdr}
\usepackage{verbatim}
\usepackage{amsfonts,longtable,mathrsfs}
\usepackage{epsfig}
\usepackage{amsfonts}
\usepackage{color}
\usepackage[numbers,sort&compress]{natbib}

\def\bc{\begin{center}}
\def\ec{\end{center}}
\newtheorem{dfn}{Definition}[section]

\begin{document}

\bc {\bf The invariance and formulas for solutions of some fifth-order difference equations
  }\ec
\medskip
\bc
D. Nyirenda %\footnote{Darlison.Nyirenda@wits.ac.za}
 and M. Folly-Gbetoula\footnote{Corresponding author: Mensah.Folly-Gbetoula@wits.ac.za} %\\Mensah.Folly-Gbetoula@wits.ac.za
 \vspace{1cm}
\\ School of Mathematics, University of the Witwatersrand, 2050, Johannesburg, South Africa.\\
 %\\e-mails: Mensah.Folly-Gbetoula@wits.ac.za, Abdul.Kara@wits.ac.za \\

\ec
\begin{abstract}
\noindent Lie group analysis of the difference equations of the form
\begin{align*}
x_{n+1} =\frac{x_{n-4}x_{n-3}}{x_{n}(a_n +b_nx_{n-4}x_{n-3}x_{n-2}x_{n-1})},
\end{align*}
where $a_n$ and $b_n$ are real sequences, is performed and non-trivial symmetries are derived. Furthermore, we find find formulas for exact solutions of the equations. This work generalizes a recent result by Elsayed [Elsayed, E.M.: {Expression and behavior of the solutions of some rational recursive sequences}. {Math. Meth. Appl. Sci.} { \bf 2016:39},5682--5694 (2016)].
\end{abstract}
\section{Introduction} \setcounter{equation}{0}
In recent years, following the work of Sophus Lie \cite{Lie} on differential equations, various researchers showed interest in symmetries. Lie investigated the group of transformations which leaves the differential equations invariant. The idea of symmetry is also connected to conservation laws and this connection between the two areas has led to greater motivation in researchers, after the work of Noether \cite{Noether}. It is known that so long as the symmetries and first integrals are related  via the invariance condition, one can implement the double reduction of the differential equations \cite{SC,MK}. The notion of using symmetries has had its extension to difference equations thanks to Maeda \cite{Maeda, Maeda2}. On symmetries in difference equations, refer to \cite{LVW,QS,hydon0,hydon1, FK, FK2, DM}. Hydon \cite{hydon0} established a symmetry based algorithm that makes solution finding possible. Despite the fact that Hydon \cite{hydon1} emphasized on lower-order difference equations, his procedure works for any order. However, for higher-order equations, computations are cumbersome as such certain assumptions are put in order to lessen the burden of computation.
\par \noindent In this paper, we are inspired by the work of Elsayed \cite{Elsayed}, who studied the following recursive sequences:
\begin{align}\label{1.0}
x_{n+1}=\frac{x_{n-3}x_{n-4}}{x_{n}(\pm 1\pm x_{n-1}x_{n-2}x_{n-3}x_{n-4})},\quad n=0,1,\dots,
\end{align}
where the initial conditions are arbitrary real numbers. Clearly, (\ref{1.0}) are special cases of a more general form
\begin{align}\label{1.1}
x_{n+1}=\frac{x_{n-3}x_{n-4}}{x_{n}(a_n +b_n x_{n-1}x_{n-2}x_{n-3}x_{n-4})},\quad n=0,1,\dots,
\end{align}
where $(a_n)$ and $(b_n)$ are real sequences. Our aim is to utilize symmetry methods to solve this more general difference equation (\ref{1.1}). Equivalently, we study the forward difference equation
\begin{align}\label{1.2}
u_{n+5}=\frac{u_nu_{n+1}}{u_{n+4}(A _n + B _n u_nu_{n+1}u_{n+2}u_{n+3})},
\end{align}
since we follow the notation of \cite{hydon0}.
%\noindent For related work, refer to \cite{ibm0}, \cite{ibm1}, \cite{ibm2} and \cite{ibm3}.
\section{Preliminaries}
This section provides background to difference equations in the context of Lie symmetry analysis.
\begin{dfn}%\cite{PO}
Let $G$ be a local group of transformations acting on a manifold $M$. A subset $\mathcal{S}\subset M$ is called $G$-invariant, and $G$ is called symmetry group of $\mathcal{S}$, if whenever $x\in \mathcal{S} $, and $g\in G$ is such that $g\cdot x$ is defined, then $g\cdot x \in \mathcal{S}$.
\end{dfn}
\begin{dfn}%\cite{PO}
Let $G$ be a connected group of transformations acting on a manifold $M$. A smooth real-valued function $\mathcal{V}: M\rightarrow \mathbb{R}$ is an invariant function for $G$ if and only if $$X(\mathcal{V})=0\qquad \text { for all } \qquad  x\in M,$$
and every infinitesimal generator $X$ of $G$.
\end{dfn}
\begin{dfn}%\cite{hydon0}
A parameterized set of point transformations,
\begin{equation}
\Gamma_{\varepsilon} :x\mapsto \hat{x}(x;\varepsilon),
\label{eq: b}
\end{equation}
where $x=x_i, $ $i=1,\dots,p$ are continuous variables, is a one-parameter local Lie group of transformations if the following conditions are satisfied:
\begin{enumerate}
\item $\Gamma_0$ is the identity map if $\hat{x}=x$ when $\varepsilon=0$
\item $\Gamma_a\Gamma_b=\Gamma_{a+b}$ for every $a$ and $b$ sufficiently close to 0
\item Each $\hat{x_i}$ can be represented as a Taylor series (in a neighborhood of $\varepsilon=0$ that is determined by $x$), and therefore
\end{enumerate}
\begin{equation}
\hat{x_i}(x:\varepsilon)=x_i+\varepsilon \xi _i(x)+O(\varepsilon ^2), i=1,...,p.
\label{eq: c}
\end{equation}
\end{dfn}
Assume  that the forward $r$th-order difference equation takes the form
\begin{align}\label{general}
%\begin{align}
u_{n+r}=&F (n, u_n, u_{n+1},\dots,  u_{n+r-1}), \quad n\in D
%\end{align}
\end{align}
for some smooth function $F$ and a regular domain $D\subset \mathbb{Z}$.
So as to compute a symmetry group of \eqref{general}, we take into consideration the group of point transformations given as
\begin{subequations}\label{Gtransfo}
\begin{align}
&\hat{n}= n,\\
&\hat{u}_n = u_n+\epsilon Q (n,u_n)+O(\epsilon ^2),\\
&\hat{u}_{n+j} = u_{n+j}+\epsilon S^j Q (n,u_n)+O(\epsilon ^2),
\end{align}
\end{subequations}
where $\epsilon$ ($\epsilon$ is sufficiently small) is the parameter, $Q=Q(n,u_n)$ is a continuous function, referred to as  characteristic and $S$ is the shift operator defined as  $S: n \mapsto n+1$. The criterion of invariance is then
\begin{align}\label{general'}
\hat{u}_{n+r}=&F (\hat{n}, \hat{u}_n, \hat{u}_{n+1},\dots,  \hat{u}_{n+r-1}),
\end{align}
which yields the linearized symmetry condition
\begin{align}\label{LSC}
 S^r Q- X F =0
\end{align}
by substituting  \eqref{Gtransfo} in \eqref{general'}. Observe that
\begin{align}\label{Ngener}
X= & Q(n,u_n)\frac{\partial}{ \partial u_n}+SQ(n,u_{n})\frac{\partial\quad }{ \partial u_{n+1}} +\dots+ S^{r-1}Q(n,u_{n})\frac{\partial\quad }{ \partial u_{n+r-1}}
\end{align}
is the corresponding $\lq$prolonged\rq \, infinitesimal of the group of transformations \eqref{Gtransfo}.
Upon knowledge of the function(s) $Q$, one is able to obtain the invariant $\mathcal{V}$  by using the canonical coordinate \cite{JV}
 \begin{align}\label{cano}
S_n= \int{\frac{du_n}{Q(n,u_n)}} .
 \end{align}
\noindent Generally, the steps involves are lengthy even though very exact and do not give room to guess work on the perfect choice of invariants. \par \noindent For more understanding on Lie analysis of differential and difference equations, see \cite{PO, hydon0}.
\section{Main results}
We are studying the equation
\begin{align}\label{un}
u_{n+5}=F =\frac{u_nu_{n+1}}{u_{n+4}(A_n +B_nu_nu_{n+1}u_{n+2}u_{n+3})}.
\end{align}
Applying condition \eqref{LSC} to \eqref{un}, we get
\begin{align}\label{a1}
 &Q(n+5, F)
 +\frac{ A_n{u_n}u_{n+1}Q\left({n+4, u_{n+4}}\right)}{{u_{n+4}}\left(A_n +B_nu_nu_{n+1}u_{n+2}u_{n+3}\right)^2 } +\nonumber\\
 &\frac{ B_n{u_n}^2{u_{n+1}}^2u_{n+2}u_{n+3}Q\left({n+4, u_{n+4}}\right)}{{u_{n+4}}^2\left(A_n +B_nu_nu_{n+1}u_{n+2}u_{n+3}\right)^2 } +%\left({n, u_{n}}\right)
%\nonumber \\&
  \frac{B_n {u_n}^2{u_{n+1}}^2u_{n+2}Q\left(n+3, u_{n+3}\right)}{u_{n+4}\left(A_n +B_nu_nu_{n+1}u_{n+2}u_{n+3}\right)^2}
  \nonumber\\&
  +\frac{B_n {u_n}^2{u_{n+1}}^2u_{n+3}Q\left(n+2, u_{n+2}\right)}{u_{n+4}\left(A_n +B_nu_nu_{n+1}u_{n+2}u_{n+3}\right)^2}  -\frac{A_nu_{n+1}Q\left(n,u_{n}\right) }{{u_{n+4}\left(A_n +B_nu_nu_{n+1}u_{n+2}u_{n+3}\right)}^{2}}=0.
\end{align}
 Eliminating $F$ is achieved by applying implicit differentiation with respect to $u_n$ (regarding $u_{n+4}$ as a function of $u_n$, $u_{n+1}$, $u_{n+2}$, $u_{n+3}$ and $u_{n+5}$) via the differential operator
 $$L=\frac{\partial\quad }{\partial {u_n}}+\frac{\partial u_{n+4}}{\partial{u_n}}\frac{\partial }{\partial u_{n+4}}
 =\frac{\partial\quad }{\partial {u_n}}-\left[\left(\frac{\partial F}{\partial{u_n}}\right)\Big/\left(\frac{\partial F}{\partial{u_{n+4}}}\right)\right]\frac{\partial\qquad }{\partial u_{n+4}}.$$
 With some simplification, we get
\begin{align}\label{a3}
& \left(A_n +B_nu_nu_{n+1}u_{n+2}u_{n+3}\right) Q'\left(n+4, u_{n+4}\right)+
B_n u_nu_{n+1}u_{n+2}Q(n+3,u_{n+3})+\nonumber\\
&
B_n u_nu_{n+1}u_{n+3}Q(n+2,u_{n+2})+
B_n u_nu_{n+2}u_{n+3}Q(n+2,u_{n+2})+\nonumber\\
&
B_n u_nu_{n+2}u_{n+3}Q(n+1,u_{n+1})
 - \left(A_n +B_nu_nu_{n+1}u_{n+2}u_{n+3}\right)Q'\left(n, u_{n}\right)+\nonumber\\
 &2 B_n u_{n+1}u_{n+2}u_{n+3} Q(n,u_{n})+
 \frac{A_n}{ u_{n}}Q\left(n, u_{n}\right) =0.
\end{align}
The symbol $'$ denotes the derivative with respect to the continuous variable. Differentiating (\ref{a3}) with respect to $u_{n}$ twice, keeping $u_{n+4}$ constant, yields
\begin{align}\label{a4'}
&-B_n u_n u_{n+1}u_{n+2}u_{n+3}Q'''(n,u_n)-A_nQ'''(n,u_n)+\frac{A_n}{u_n}
Q''\left(n,{u_{n}}\right)\nonumber\\
 &-\frac{2A_n}{{u_n}^2}Q'\left(n,u_{n}\right)+ \frac{2A_n}{{u_n}^3}Q\left(n,u_{n}\right)=0.
\end{align}
The characteristic in \eqref{a4'} is a function of $u_{n}$ only and thus we split \eqref{a4'} to get the system
%$$
\begin{subequations}\label{a4}
\begin{align}
%\begin{cases}
%\begin{align}\label{a4}
 &1:Q'''(n,u_n)-\frac{1}{u_n}
Q''\left(n,{u_{n}}\right)+\frac{2}{{u_n}^2}Q'\left(n,u_{n}\right)
 - \frac{2}{{u_n}^3}Q\left(n,u_{n}\right)=0\\&\nonumber\\
& u_{n+1}u_{n+2}u_{n+3}:Q'''(n,u_n)=0.
%\end{align}
% \end{cases}
 \end{align}
%$$
\end{subequations}
One obtains the solution to \eqref{a4} as
\begin{align}\label{a6}
\begin{split}
Q\left(n,{u_{n}}\right) = \alpha_ n  {u_n}^2 +\beta _n  {u_n}
\end{split}
\end{align}
for some arbitrary functions $\alpha _n$ and $\beta _n$ of $n$. Substituting  \eqref{a6} and its  shifts  in \eqref{a1}, and making a replacement of the expression of $u_{n+5}$ given in \eqref{un} in the resulting equation leads to
\begin{align}\label{a7'}
& B_n u_nu_{n+1}u_{n+2}{u_{n+3}}^2u_{n+4}\alpha_{n+3}+
B_n u_nu_{n+1}u_{n+2}{u_{n+3}}{u_{n+4}}^2\alpha_{n+4}+\nonumber\\
&B_n u_nu_{n+1}{u_{n+2}}^2{u_{n+3}}u_{n+4}\alpha_{n+2}+
B_n u_nu_{n+1}u_{n+2}{u_{n+3}}u_{n+4}(\beta_{n+2}+\beta_{n+3}+\nonumber\\&
\beta_{n+4}+\beta_{n+5})
+A_n {u_{n+4}}^2\alpha_{n+4}
-A_n {u_{n}}{u_{n+4}}\alpha_{n}
-A_n {u_{n+1}}{u_{n+4}}\alpha_{n+1}-\nonumber\\
&
A_nu_{n+4}\left(  \beta_{n} +\beta_{n+1}-\beta_{n+4}-\beta_{n+5} \right)+u_nu_{n+1}\alpha_{n+5}=0.
\end{align}
Now equate coefficients of all powers of shifts of $u_n$ to zero, i.e.,
\begin{align}\label{a7}
&u_nu_{n+1}u_{n+2}{u_{n+3}}^2u_{n+4}&:&\quad \alpha_{n+3}=0\\
&u_nu_{n+1}u_{n+2}{u_{n+3}}{u_{n+4}}^2&:&\quad\alpha_{n+4}=0\\
&u_nu_{n+1}{u_{n+2}}^2{u_{n+3}}u_{n+4}&:&\quad\alpha_{n+2}=0\\
& u_nu_{n+1}u_{n+2}{u_{n+3}}u_{n+4}&:&\quad(\beta_{n+2}+\beta_{n+3}+
\beta_{n+4}+\beta_{n+5})=0\\
&{u_{n+4}}^2&:&\quad\alpha_{n+4}=0\\
& {u_{n}}{u_{n+4}}&:&\quad\alpha_{n}=0\\
& {u_{n+1}}{u_{n+4}}&:&\quad\alpha_{n+1}=0\\
&u_{n+4}&:&\quad\left(  \beta_{n} +\beta_{n+1}-\beta_{n+4}-\beta_{n+5} \right)=0\\
&u_nu_{n+1}&:&\quad\alpha_{n+5}=0.
\end{align}
So the system above is reduced to
\begin{align}\label{rel}
&\alpha _n=0,  \\
&\beta_{n} + \beta _{n+1}+ \beta _{n+2}+ \beta _{n+3} =0\label{rel}.
\end{align}
The three independent solutions of the linear third-order difference equation \eqref{rel} are given by
\begin{align}
(-1)^n, \quad \beta ^n  \quad \text{ and  } \quad \bar{\beta}^n
\end{align}
where $\beta = \exp \{i\pi/2\}$ and $\bar{\beta}$ denotes its complex conjugate.
The characteristics are then given by
\begin{align}\label{a9}
Q_1(n,u_n)=(-1) ^n u_n, \quad Q_2(n,u_n)=\beta ^n u_n \quad \text{ and  }\quad  Q_3(n,u_n)=\bar{\beta}^n u_n,
\end{align}
and therefore,  the symmetry operators admitted by \eqref{un} are given by
\begin{align}\label{gener}
X_1= &\sum\limits_{j=0}^{4}(-1) ^{n+j} u_{n+j} \frac{\partial}{\partial _{ u_{n+j}}},\; X_2= \sum\limits_{j=0}^{4}\beta ^{n+j} u_{n+j} \frac{\partial}{\partial _{ u_{n+j}}},\nonumber\\ X_3= &  \sum\limits_{j=0}^{4}\bar{\beta} ^{n+j} u_{n+j} \frac{\partial}{\partial _{ u_{n+j}}}.
\end{align}
One can choose any one of the characteristics to write  the canonical coordinate. We select $Q_2$. Thus
\begin{align}\label{cano}
S_n =\int\frac{du_n}{Q_2(n,u_n)}=\int\frac{du_n}{\beta ^nu_n}=\frac{1}{\beta ^n}\ln|u_n|
\end{align}
and we use relation \eqref{rel} to derive the invariant function $\tilde{V}_n$ as follows:
\begin{align}\label{tilde}
\tilde{V}_n = S_n\beta^ n +S_{n+1}\beta^{n+1}+S_{n+2}\beta^{n+2}+ S_{n+3}\beta^{n+3}.
\end{align}
Actually,
\begin{align}
X_1 (\tilde{V}_n) = (-1)^ n+ (-1)^{n+1}+ (-1)^{n+2}+ (-1)^{n+3}=0,
\end{align}
\begin{align}
X_2 (\tilde{V}_n) = \beta^ n+ \beta^{n+1}+ \beta^{n+2}+ \beta^{n+3}=0
\end{align}
and
\begin{align}
X_3 (\tilde{V}_n) = \bar{\beta}^ n  + \bar{\beta}^{n+1} + \bar{\beta}^{n+2}+ \bar{\beta}^{n+3}=0.
\end{align}
For the sake of simplicity, we utilize
\begin{align}\label{vn}
|{V}_n| =\exp\{  -\tilde{V}_n \}
\end{align}
instead. In other words, $V_n =\pm 1/(u_nu_{n+1}u_{n+2}u_{n+3})$.
One can show  via \eqref{un} and \eqref{vn} that
\begin{equation}\label{vn1}
V_{n+2}={A_n V_n}\pm B_n.%\qquad \text{with}\quad v_0=u_0u_1\neq -1/c.
\end{equation}
By utilizing the plus sign (one is allowed to choose), the solution of \eqref{vn1} can be presented in closed form as follows:
\begin{align}\label{solvn}
V_{2n+j}\quad=& V_j \left(   \prod\limits_{k_1=0}^{n-1}A_{2k_1+j}\right) +\sum\limits _{l=0}^{n-1} \left(  B_{2l+j}\prod\limits _{k_2=l+1}^{n-1}A_{2k_2}\right), \quad j=0,1.
\end{align}
From the above equation, obtaining the solution of \eqref{un} is easier. We first use \eqref{cano} to get
\begin{align}
\vert u_{n}\vert  =& \exp\left(\beta _n S_n\right).
\end{align}
Secondly, we use \eqref{tilde} to get
\begin{align}
\vert u_{n}\vert
=&\exp\Bigg[\beta^{n}c_1 + \bar{\beta}^{n}c_2 +(-1)^nc_3- \left(\frac{1}{4}-\frac{i}{4}\right)\sum\limits_{k_1 = 0}^{n - 1}{\beta}^{n}\bar{\beta}^{k_1}|{\tilde{V}}_{k_1}|\nonumber\\
            &- \left(\frac{1}{4}+\frac{i}{4}\right)\sum\limits_{k_2 = 0}^{n - 1} \bar{\beta}^{n}{\beta}^{k_2}|\tilde{V}_{k_2}|-
           \frac{1}{2}\sum\limits_{k_3=0}^{n-1}
           (-1)^{n-k_3} |{\tilde{V}}_{k_3}|\Bigg].
\end{align}
Finally, invoking \eqref{vn} yields
\begin{align}
\vert u_{n} \vert
           =& \exp\Bigg[\beta^{n}c_1 + \bar{\beta}^{n}c_2 +(-1)^nc_3+ \left(\frac{1}{4}-\frac{i}{4}\right)\sum\limits_{k_1 = 0}^{n - 1}{\beta}^{n}\bar{\beta}^{k_1}\ln|V_{k_1}|\nonumber\\
            &+ \left(\frac{1}{4}+\frac{i}{4}\right)\sum\limits_{k_2 = 0}^{n - 1} \bar{\beta}^{n}{\beta}^{k_2}\ln| V_{k_2}|+
           \frac{1}{2}\sum\limits_{k_3=0}^{n-1}
           (-1)^{n-k_3}\ln |V_{k_3}\Bigg]\nonumber\\
           = & \exp\Bigg( H_{n} + \frac{1}{2}\sum\limits_{k = 0}^{n - 1}\left[ \sqrt{2}\cos\left(\frac{\pi(2k - 2n + 1)}{4}\right) + (-1)^{k - n}\right]\ln|V_{k}|\Bigg),\label{unsol}
\end{align}
where $H_{n} = \beta^{n}c_1 + \bar{\beta}^{n}c_2 +(-1)^nc_3$.
Replacing $n$ with $4n + j$ for $j = 0,1,2,3$ yields
\begin{equation}\label{main}
\vert u_{4n + j}\vert  = \exp\Bigg[ H_{j} + \frac{1}{2}\sum\limits_{k = 0}^{n - 1}\left( \sqrt{2}\cos\left(\frac{\pi(2k - 2j + 1)}{4}\right) + (-1)^{k - j}\right)\ln|V_{k}|\Bigg].
\end{equation}
Set $j = 0$ in \eqref{main} to get
$$\vert   u_{4n} \vert = \exp(H_{0})\prod\limits_{s = 0}^{n - 1}\vert \frac{V_{4s}}{V_{4s + 1}}\vert.$$
\noindent But substituting $n = 0$ in \eqref{unsol} leads to $\vert u_0 \vert = \exp(H_{0})$. Furthermore, using \eqref{un} and \eqref{vn}, it can be shown that there is no need of absolute values. Hence
\begin{align*}\label{eqn0}
u_{4n} & = u_0\prod\limits_{s = 0}^{n - 1}\frac{V_{4s}}{V_{4s + 1}}\\
                   & = u_0\prod\limits_{s = 0}^{n - 1}\frac{ V_0 \left(   \prod\limits_{k_1=0}^{2s-1}A_{2k_1}\right) +\sum\limits _{l=0}^{2s-1} \left(  B_{2l}\prod\limits _{k_2=l+1}^{2s-1}A_{2k_2}\right)}{V_1 \left(   \prod\limits_{k_1=0}^{2s-1}A_{2k_1+1}\right) +\sum\limits_{l=0}^{2s-1} \left(  B_{2l+1}\prod\limits_{k_2=l+1}^{2s-1}A_{2k_2}\right)}\\
                   & = \frac{u_4^{n}}{u_0^{n - 1}}\prod\limits_{s = 0}^{n - 1}\frac{\left(   \prod\limits_{k_1=0}^{2s-1}A_{2k_1}\right) + u_0u_1u_2u_3\sum\limits_{l=0}^{2s-1} \left(  B_{2l}\prod\limits_{k_2=l+1}^{2s-1}A_{2k_2}\right)}{\left(   \prod\limits_{k_1=0}^{2s-1}A_{2k_1+1}\right) + u_1u_2u_3u_4\sum\limits_{l=0}^{2s-1} \left(  B_{2l+1}\prod\limits_{k_2=l+1}^{2s-1}A_{2k_2 + 1}\right)}.
\end{align*}
For $j = 1$, we find that
$$ u_{4n + 1} = u_1\prod_{s = 0}^{n - 1}\frac{V_{4s + 1}}{V_{4s + 2}}$$ so that
\begin{align*}
u_{4n + 1} & = u_1\prod\limits_{s = 0}^{n - 1}\frac{ V_1 \left(\prod\limits_{k_1=0}^{2s-1}A_{2k_1+1}\right) +\sum\limits_{l=0}^{n-1} \left(  B_{2l+1}\prod\limits_{k_2=l+1}^{2s-1}A_{2k_2 + 1}\right) }{ V_0 \left( \prod\limits_{k_1=0}^{2s}A_{2k_1}\right) +\sum\limits_{l=0}^{2s} \left(  B_{2l}\prod\limits_{k_2=l+1}^{2s}A_{2k_2}\right)}\\
           & = \frac{u_0^{n}u_1}{u_4^{n}}\prod\limits_{s = 0}^{n - 1}\frac{\left(   \prod\limits _{k_1=0}^{2s-1}A_{2k_1+1}\right) + u_1u_2u_3u_4\sum\limits _{l=0}^{2s-1} \left(  B_{2l+1}\prod\limits _{k_2=l+1}^{2s-1}A_{2k_2 + 1}\right) }{ \left(   \prod\limits _{k_1=0}^{2s}A_{2k_1}\right) + u_0u_1u_2u_3\sum\limits _{l=0}^{2s} \left(B_{2l}\prod\limits _{k_2=l+1}^{2s}A_{2k_2}\right)}.
\end{align*}
For $j = 2$, we have
$$ u_{4n + 2} = u_2\prod\limits_{s = 0}^{n - 1}\frac{V_{4s + 2}}{V_{4s + 3}}$$ which evaluates to
\begin{align*}
u_{4n + 2} & = u_2\prod\limits_{s = 0}^{n - 1}\frac{ V_0 \left( \prod\limits_{k_1=0}^{2s}A_{2k_1}\right) +\sum\limits _{l=0}^{2s} \left(  B_{2l}\prod\limits _{k_2=l+1}^{2s}A_{2k_2}\right) }{V_1 \left(   \prod\limits _{k_1=0}^{2s}A_{2k_1+1}\right) +\sum\limits _{l=0}^{2s} \left(  B_{2l+1}\prod\limits_{k_2=l+1}^{2s}A_{2k_2 + 1}\right)}\\
          & = \frac{u_4^{n}u_2}{u_0^{n}}\prod\limits_{s = 0}^{n - 1}\frac{\left(   \prod\limits _{k_1=0}^{2s}A_{2k_1}\right) + u_0u_1u_2u_3\sum\limits_{l=0}^{2s} \left(  B_{2l}\prod\limits_{k_2=l+1}^{2s}A_{2k_2}\right) }{\left(   \prod\limits_{k_1=0}^{2s}A_{2k_1+1}\right) + u_1u_2u_3u_4\sum\limits_{l=0}^{2s} \left(  B_{2l+1}\prod\limits_{k_2=l+1}^{2s}A_{2k_2 + 1}\right)}.
\end{align*}
Finally, for $j = 3$, we obtain
$$ u_{4n + 3} = u_3\prod\limits_{s = 0}^{n - 1}\frac{V_{4s + 3}}{V_{4s + 4}}$$ so that
\begin{align*}
u_{4n + 3} & = u_3\prod\limits_{s = 0}^{n - 1}\frac{ V_1 \left(\prod\limits _{k_1=0}^{2s}A_{2k_1+1}\right) +\sum\limits_{l=0}^{2s} \left(  B_{2l+1}\prod\limits_{k_2=l+1}^{2s}A_{2k_2 + 1}\right) }{V_0 \left(   \prod\limits_{k_1=0}^{2s + 1}A_{2k_1}\right) +\sum\limits_{l=0}^{2s + 1} \left(  B_{2l}\prod\limits_{k_2=l+1}^{2s + 1}A_{2k_2}\right)}\\
          & = \frac{u_0^{n}u_3}{u_4^{n}}\prod\limits_{s = 0}^{n - 1}\frac{\left(   \prod\limits_{k_1=0}^{2s}A_{2k_1+1}\right) + u_1u_2u_3u_4\sum _{l=0}^{2s} \left(  B_{2l+1}\prod\limits_{k_2=l+1}^{2s}A_{2k_2 + 1}\right) }{\left(   \prod\limits_{k_1=0}^{2s + 1}A_{2k_1}\right) + u_0u_1u_2u_3\sum _{l=0}^{2s + 1} \left(  B_{2l}\prod\limits_{k_2=l+1}^{2s + 1}A_{2k_2}\right)}.
\end{align*}
Hence the solution to \eqref{1.1} is given by
\begin{equation*}
x_{4n - 4} = \frac{x_{0}^{n}}{x_{-4}^{n - 1}}\prod\limits_{s = 0}^{n - 1}\frac{\left(   \prod\limits_{k_1=0}^{2s-1}a_{2k_1}\right) + x_{-4}x_{-3}x_{-2}x_{-1}\sum\limits_{l=0}^{2s-1} \left(  b_{2l}\prod\limits_{k_2=l+1}^{2s-1}a_{2k_2}\right)}{\left(   \prod\limits_{k_1=0}^{2s-1}a_{2k_1+1}\right) + x_{-3}x_{-2}x_{-1}x_0\sum _{l=0}^{2s-1} \left(  b_{2l+1}\prod\limits_{k_2=l+1}^{2s-1}a_{2k_2 + 1}\right)},
\end{equation*}
which can be rearranged as
\begin{equation*}
x_{4n} = \frac{x_{0}^{n + 1}}{x_{-4}^{n}}\prod\limits_{s = 0}^{n}\frac{\left(   \prod\limits_{k_1=0}^{2s-1}a_{2k_1}\right) + x_{-4}x_{-3}x_{-2}x_{-1}\sum\limits_{l=0}^{2s-1} \left(  b_{2l}\prod\limits_{k_2=l+1}^{2s-1}a_{2k_2}\right)}{\left(   \prod\limits_{k_1=0}^{2s-1}a_{2k_1+1}\right) + x_{-3}x_{-2}x_{-1}x_0\sum _{l=0}^{2s-1} \left(  b_{2l+1}\prod\limits_{k_2=l+1}^{2s-1}a_{2k_2 + 1}\right)}.
\end{equation*}
The term $s = 0$ in the product (indexed by $s$) is equal to 1 using the facts that $\sum\limits_{i = 0}^{-1}a_i = 0$ and $\prod\limits_{j = 0}^{-1}a_j = 1$. As a result, we can still rewrite the solution as
\begin{subequations}\label{xnsol}
\begin{equation}\label{sol0}
x_{4n} = \frac{x_{0}^{n + 1}}{x_{-4}^{n}}\prod\limits_{s = 0}^{n - 1}\frac{\left(   \prod\limits_{k_1=0}^{2s + 1}a_{2k_1}\right) + x_{-4}x_{-3}x_{-2}x_{-1}\sum\limits_{l=0}^{2s + 1} \left(  b_{2l}\prod\limits_{k_2=l+1}^{2s + 1}a_{2k_2}\right)}{\left(   \prod\limits_{k_1=0}^{2s + 1}a_{2k_1+1}\right) + x_{-3}x_{-2}x_{-1}x_0\sum _{l=0}^{2s + 1} \left(  b_{2l+1}\prod\limits_{k_2=l+1}^{2s + 1}a_{2k_2 + 1}\right)}.
\end{equation}
Furthermore, observe that
\begin{equation}\label{sol1}
x_{4n - 3} = \frac{x_{-4}^{n}x_{-3}}{x_{0}^{n}}\prod\limits_{s = 0}^{n - 1}\frac{\left(   \prod\limits_{k_1=0}^{2s-1}a_{2k_1+1}\right) + x_{-3}x_{-2}x_{-1}x_0 \sum\limits_{l=0}^{2s-1} \left(  b_{2l+1}\prod\limits_{k_2=l+1}^{2s-1}a_{2k_2 + 1}\right) }{ \left(   \prod\limits_{k_1=0}^{2s}a_{2k_1}\right) + x_{-4}x_{-3}x_{-2}x_{-1}  \sum\limits_{l=0}^{2s} \left(b_{2l}\prod\limits_{k_2=l+1}^{2s}a_{2k_2}\right)},
\end{equation}
\begin{equation}\label{sol2}
x_{4n  - 2} = \frac{x_{0}^{n}x_{-2}}{x_{-4}^{n}}\prod\limits_{s = 0}^{n - 1}\frac{\left(   \prod\limits_{k_1=0}^{2s}a_{2k_1}\right) + x_{-4}x_{-3}x_{-2}x_{-1}\sum\limits_{l=0}^{2s} \left(  b_{2l}\prod\limits_{k_2=l+1}^{2s}a_{2k_2}\right) }{\left(   \prod\limits_{k_1=0}^{2s}a_{2k_1+1}\right) + x_{-3}x_{-2}x_{-1}x_0\sum\limits_{l=0}^{2s} \left(  b_{2l+1}\prod\limits_{k_2=l+1}^{2s}a_{2k_2 + 1}\right)}
\end{equation}
and
\begin{equation}\label{sol3}
x_{4n - 1} = \frac{x_{-4}^{n}x_{-1}}{x_{0}^{n}}\prod\limits_{s = 0}^{n - 1}\frac{\left(   \prod\limits_{k_1=0}^{2s}a_{2k_1+1}\right) + x_{-3}x_{-2}x_{-1}x_0\sum\limits_{l=0}^{2s} \left(  b_{2l+1}\prod\limits_{k_2=l+1}^{2s}a_{2k_2 + 1}\right) }{\left(   \prod\limits _{k_1=0}^{2s + 1}a_{2k_1}\right) +  x_{-4}x_{-3}x_{-2}x_{-1}\sum\limits_{l=0}^{2s + 1} \left(  b_{2l}\prod\limits_{k_2=l+1}^{2s + 1}a_{2k_2}\right)}
\end{equation}
\end{subequations}
as long as any of the denominators does not vanish.\\
\noindent In the following sections, we look at some special cases.
\section{The case $a_n$, $b_n$ are 1-periodic}
In this case $a_n = a$ and $b_n = b$ where $a, b \in \mathbb{R}$.
\subsection{The case $a \neq 1$}
From \eqref{xnsol}, the solution is given by
\begin{subequations}\label{neq1}
\begin{equation}\label{ssol0}
x_{4n} = \frac{x_{0}^{n + 1}}{x_{-4}^{n}}\prod\limits_{s = 0}^{n - 1}\frac{ a^{2s + 2} + bx_{-4}x_{-3}x_{-2}x_{-1}\frac{1 - a^{2s + 2}}{1 - a}}{a^{2s + 2} + bx_{-3}x_{-2}x_{-1}x_0\frac{1 - a^{2s + 2}}{1 - a}},
\end{equation}
\begin{equation}\label{ssol1}
x_{4n - 3} = \frac{x_{-4}^{n}x_{-3}}{x_{0}^{n}}\prod\limits_{s = 0}^{n - 1}\frac{  a^{2s} + bx_{-3}x_{-2}x_{-1}x_0\frac{1 - a^{2s}}{1 - a}}{ a^{2s + 1} + bx_{-4}x_{-3}x_{-2}x_{-1}  \frac{1 - a^{2s + 1}}{1 - a}},
\end{equation}
\begin{equation}\label{ssol2}
x_{4n  - 2} = \frac{x_{0}^{n}x_{-2}}{x_{-4}^{n}}\prod\limits_{s = 0}^{n - 1}\frac{  a^{2s + 1} + bx_{-4}x_{-3}x_{-2}x_{-1}\frac{1 - a^{2s + 1}}{1 - a} }{ a^{2s + 1} + bx_{-3}x_{-2}x_{-1}x_0 \frac{1 - a^{2s + 1}}{1 - a}}
\end{equation}
and
\begin{equation}\label{ssol3}
x_{4n - 1} = \frac{x_{-4}^{n}x_{-1}}{x_{0}^{n}}\prod\limits_{s = 0}^{n - 1}\frac{ a^{2s + 1} + bx_{-3}x_{-2}x_{-1}x_0 \frac{1 - a^{2s + 1}}{1 - a}   }{ a^{2s + 2} +  bx_{-4}x_{-3}x_{-2}x_{-1} \frac{1 - a^{2s + 2}}{1 - a}}
\end{equation}
\end{subequations}
where $x_{-4}, x_0 \neq 0$ and for all $(i,s) \in \{0,1\} \times \{0,1,2,3,\ldots, n - 1\}$, $$(1 - a)a^{2s + i} + (1 - a^{2s + i})bx_{-3}x_{-2}x_{-1}x_{0} \neq 0$$
 and
$$(1 - a)a^{2s + 1 + i} + (1 - a^{2s + 1 + i})bx_{-4}x_{-3}x_{-2}x_{-1}\neq 0.$$
\subsubsection{ The case $a = -1$}
In this case, the solution which for $b = \pm 1$ appears in \cite{Elsayed} (see Theorems 3 and 8), is given by
\begin{subequations}\label{n-1}
\begin{equation}\label{sl0}
x_{4n} = \frac{x_{0}^{n + 1}}{x_{-4}^{n}},
\end{equation}
\begin{equation}\label{sl1}
x_{4n - 3} = \frac{x_{-4}^{n}x_{-3}}{x_{0}^{n}}(-1 + bx_{-4}x_{-3}x_{-2}x_{-1})^{-n},
\end{equation}
\begin{equation}\label{sl2}
x_{4n  - 2} = \frac{x_{0}^{n}x_{-2}}{x_{-4}^{n}}\left(\frac{  -1 + bx_{-4}x_{-3}x_{-2}x_{-1} }{ -1 + bx_{-3}x_{-2}x_{-1}x_0 }\right)^{n}
\end{equation}
and
\begin{equation}\label{sl3}
x_{4n - 1} = \frac{x_{-4}^{n}x_{-1}}{x_{0}^{n}}(-1 + bx_{-3}x_{-2}x_{-1}x_0)^{n}
\end{equation}
\end{subequations}
where $x_{-4}, x_0 \neq 0$, $bx_{-4}x_{-3}x_{-2}x_{-1} \neq 1$ and $bx_{-3}x_{-2}x_{-1}x_{0} \neq 1$.
\subsection{The case $a = 1$}
From \eqref{xnsol}, the solution, which for $b = \pm 1$ appears in \cite{Elsayed} (see Theorems 1 and 6), is given by
\begin{subequations}\label{n1}
\begin{equation}\label{l0}
x_{4n} = \frac{x_{0}^{n + 1}}{x_{-4}^{n}}\prod\limits_{s = 0}^{n - 1}\frac{ 1 + (2s + 2)bx_{-4}x_{-3}x_{-2}x_{-1}}{1 + (2s + 2)bx_{-3}x_{-2}x_{-1}x_0},
\end{equation}
\begin{equation}\label{l1}
x_{4n - 3} = \frac{x_{-4}^{n}x_{-3}}{x_{0}^{n}}\prod\limits_{s = 0}^{n - 1}\frac{ 1 + 2sbx_{-3}x_{-2}x_{-1}x_0 }{1 + (2s + 1)bx_{-4}x_{-3}x_{-2}x_{-1}},
\end{equation}
\begin{equation}\label{l2}
x_{4n  - 2} = \frac{x_{0}^{n}x_{-2}}{x_{-4}^{n}}\prod\limits_{s = 0}^{n - 1}\frac{1 + (2s + 1)bx_{-4}x_{-3}x_{-2}x_{-1}}{ 1 + (2s + 1)bx_{-3}x_{-2}x_{-1}x_0}
\end{equation}
and
\begin{equation}\label{l3}
x_{4n - 1} = \frac{x_{-4}^{n}x_{-1}}{x_{0}^{n}}\prod\limits_{s = 0}^{n - 1}\frac{ 1 + (2s + 1)bx_{-3}x_{-2}x_{-1}x_0}{ 1 +  (2s + 2)bx_{-4}x_{-3}x_{-2}x_{-1}}
\end{equation}
\end{subequations}
where $x_{-4}, x_0 \neq 0$, $2jbx_{-4}x_{-3}x_{-2}x_{-1} \neq -1$, $(2j - 1)bx_{-4}x_{-3}x_{-2}x_{-1} \neq -1$, $2jbx_{-3}x_{-2}x_{-1}x_{0} \neq -1$ and
$(2j - 1)bx_{-3}x_{-2}x_{-1}x_{0} \neq -1$ for all $j = 1,2,3, \ldots, n$.
\section{The case $a_n, b_n$ are 2-periodic}
We assume that $\{ a_n\}_{n = 0}^{\infty} = a_0, a_1, a_0, a_1, \ldots$ and $\{ b_n\}_{n = 0}^{\infty} = b_0, b_1, b_0, b_1, \ldots$.
Then, from \eqref{xnsol}, we have
\begin{subequations}\label{ntwo}
\begin{equation}\label{sol00}
x_{4n} = \frac{x_{0}^{n + 1}}{x_{-4}^{n}}\prod\limits_{s = 0}^{n - 1}\frac{a_{0}^{2s + 2} + b_0x_{-4}x_{-3}x_{-2}x_{-1}\sum\limits_{l=0}^{2s + 1}a_0^{l}}{ a_{1}^{2s + 2} + b_1x_{-3}x_{-2}x_{-1}x_0\sum\limits_{l=0}^{2s + 1}a_1^{l}},
\end{equation}
\begin{equation}\label{sol10}
x_{4n - 3} = \frac{x_{-4}^{n}x_{-3}}{x_{0}^{n}}\prod\limits_{s = 0}^{n - 1}\frac{a_{1}^{2s} + b_1x_{-3}x_{-2}x_{-1}x_0 \sum\limits_{l=0}^{2s-1}a_1^{l} }{ a_{0}^{2s + 1} + b_0x_{-4}x_{-3}x_{-2}x_{-1}  \sum\limits _{l=0}^{2s}a_0^{l}},
\end{equation}
\begin{equation}\label{sol20}
x_{4n  - 2} = \frac{x_{0}^{n}x_{-2}}{x_{-4}^{n}}\prod\limits_{s = 0}^{n - 1}\frac{ a_{0}^{2s + 1} + b_0x_{-4}x_{-3}x_{-2}x_{-1}\sum\limits_{l=0}^{2s}a_0^{l}}{a_{1}^{2s + 1} + b_1x_{-3}x_{-2}x_{-1}x_0\sum\limits_{l=0}^{2s}a_1^{l}}
\end{equation}
and
\begin{equation}\label{sol30}
x_{4n - 1} = \frac{x_{-4}^{n}x_{-1}}{x_{0}^{n}}\prod\limits_{s = 0}^{n - 1}\frac{a_{1}^{2s + 1} + b_1x_{-3}x_{-2}x_{-1}x_0\sum\limits_{l=0}^{2s}a_1^{l} }{ a_{0}^{2s + 2} +  b_0x_{-4}x_{-3}x_{-2}x_{-1}\sum\limits_{l=0}^{2s + 1}a_0^{l}}
\end{equation}
\end{subequations}
as long as $x_{-4}, x_0 \neq 0$ and for all $(i,s)\in \{0,1\} \times \{0,1,2,\ldots, n - 1\}$,
$a_0^{2s + 1 + i} + b_0x_{-4}x_{-3}x_{-2}x_{-1}\sum\limits_{l = 0}^{2s + i}a_0^{l} \neq 0$  and
$a_1^{2s + 1 + i} + b_1x_{-3}x_{-2}x_{-1}x_0\sum\limits_{l = 0}^{2s + i}a_1^{l} \neq 0$.
\subsection{The case $a_0 = 1$ and $a_1 = -1$}
The solution is given by
\begin{subequations}\label{ntwoo}
\begin{equation}\label{s00}
x_{4n} = \frac{x_{0}^{n + 1}}{x_{-4}^{n}}\prod\limits_{s = 0}^{n - 1}(1 + (2s + 2)b_0x_{-4}x_{-3}x_{-2}x_{-1}),
\end{equation}
\begin{equation}\label{s10}
x_{4n - 3} = \frac{x_{-4}^{n}x_{-3}}{x_{0}^{n}}\prod\limits_{s = 0}^{n - 1}\frac{1}{1 + (2s + 1)b_0x_{-4}x_{-3}x_{-2}x_{-1}},
\end{equation}
\begin{equation}\label{s20}
x_{4n  - 2} = x_{-2}\left(\frac{x_{0}}{x_{-4}(-1 + b_1x_{-3}x_{-2}x_{-1}x_0 )}\right)^{n}\prod\limits_{s = 0}^{n - 1}(1 + (2s + 1)b_0x_{-4}x_{-3}x_{-2}x_{-1}),
\end{equation}
and
\begin{equation}\label{s30}
x_{4n - 1} = x_{-1}\left(\frac{x_{-4}(-1 + b_1x_{-3}x_{-2}x_{-1}x_0)}{x_{0}}\right)^{n}\prod\limits_{s = 0}^{n - 1}\frac{1}{1 +  (2s + 2)b_0x_{-4}x_{-3}x_{-2}x_{-1}}
\end{equation}
\end{subequations}
where $x_{-4}, x_0 \neq 0$, $b_1x_{-3}x_{-2}x_{-1}x_{0} \neq 1$ and $jb_0x_{-4}x_{-3}x_{-2}x_{-1} \neq -1$ for all $j = 1,2,3,\ldots, 2n$.
\subsection{The case $a_0 = -1$ and $a_1 = 1$}
In this case, we obtain
\begin{subequations}
\begin{equation}\label{eq00}
x_{4n} = \frac{x_{0}^{n + 1}}{x_{-4}^{n}}\prod\limits_{s = 0}^{n - 1}\frac{1}{ 1 + (2s + 2)b_1x_{-3}x_{-2}x_{-1}x_0},
\end{equation}
\begin{equation}\label{eq10}
x_{4n - 3} = x_{-3}\left(\frac{x_{-4}}{x_{0}(-1 + b_0x_{-4}x_{-3}x_{-2}x_{-1})}\right)^{n}\prod\limits_{s = 0}^{n - 1}(1 + 2sb_1x_{-3}x_{-2}x_{-1}x_0),
\end{equation}
\begin{equation}\label{eq20}
x_{4n  - 2} = x_{-2}\left(\frac{x_{0}(-1 + b_0x_{-4}x_{-3}x_{-2}x_{-1})}{x_{-4}}\right)^{n}\prod\limits_{s = 0}^{n - 1}\frac{1}{1 + (2s + 1)b_1x_{-3}x_{-2}x_{-1}x_0}
\end{equation}
and
\begin{equation}\label{eq30}
x_{4n - 1} = \frac{x_{-4}^{n}x_{-1}}{x_{0}^{n}}\prod\limits_{s = 0}^{n - 1}(1  + (2s + 1)b_1x_{-3}x_{-2}x_{-1}x_0)
\end{equation}
\end{subequations}
where $x_{-4}, x_0 \neq 0$, $b_0x_{-4}x_{-3}x_{-2}x_{-1} \neq 1$ and $jb_1x_{-3}x_{-2}x_{-1}x_0 \neq -1$ for all $j = 1,2,3,\ldots, 2n$.
\section{Conclusion}
Our work in this paper was twofold. First, we found non-trivial Lie symmetry generators of the difference equations \eqref{1.1}. Second, we derived explicit formulas for solutions of difference equations in question. Consequently, this generalised what Elsayed found in \cite{Elsayed} where the values of $a_n$ and $b_{n}$ were only confined to $\pm 1$. We showed that in those particular cases, our results yielded Elsayed's results.
%%%%%%%%%%%%%%%%%%%%%%%%%%%%%%

\end{document}